\documentclass[11pt]{article}
\usepackage{amsmath}
\usepackage{amsfonts}
\usepackage{amsthm}
\usepackage{amssymb}
\usepackage{latexsym}
\usepackage{euscript}
\usepackage{xypic}
\usepackage{verbatim}
\usepackage{mathrsfs}
\usepackage{aeguill}

\newtheorem{theorem}{Theorem}[section]
\newtheorem{lemma}[theorem]{Lemma}
\newtheorem{df}[theorem]{Definition}
\newtheorem{cor}[theorem]{Corollary}

\def\vpi{g}

\def\vU{U_t}

\def\Sp{\mathrm{Spf}}

\def\q{\mathbf{q}}
\def\B{\mathcal{B}}
\def\cyc{\chi}

\def\E{\mathcal{E}}
\def\wt{\widetilde}
\def\lll{\langle \kern-0.05em{\langle}}
\def\rrr{\rangle \kern-0.05em{\rangle}}
\def\Gal{\mathrm{Gal}}
\def\Q{\mathbf{Q}}
\def\Qbar{\overline{\Q}}
\def\Z{\mathbf{Z}}

\def\W{\mathcal W}

\def\eps{\epsilon}
\def\C{\mathbf{C}}

\def\F{\mathbf{F}}
\def\sF{\mathcal{F}}
\def\sheaF{\sF}
\def\sheaG{\mathcal{G}}

\def\Fbar{\overline{\F}}
\def\Zbar{\overline{\Z}}
\def\VA{V\kern-0.15em{A}}
\begin{document}

\title{The Eigencurve is Proper at Integral Weights}
\author{Frank Calegari\footnote{Supported in part by the
American Institute of Mathematics.}}
\maketitle

\section{Introduction}

The eigencurve $\E$ is a rigid analytic space parameterizing
overconvergent and therefore classical modular eigenforms of finite slope.
Since Coleman and Mazur's original work~\cite{eigencurve}, there
have been numerous 
generalizations~\cite{eigenvarieties,chenevier,urban},
as well as
alternative constructions using modular symbols~\cite{ashstevens}
and $p$-adic representation theory~\cite{emerton}.
In spite of these advances, several elementary questions about the
geometry of $\E$ remain. One such 
question was raised by
Coleman and Mazur:
does there exist a $p$-adic family
of finite slope overconvergent eigenforms over
a punctured disk, and converging, at the puncture,
to an overconvergent eigenform of infinite slope?
Another way of phrasing this question is to ask whether
the projection $\pi: \E \rightarrow \W$ satisfies
the valuative criterion for properness\footnote{The curve
 $\E$ has infinite degree over weight space
$\W$,
and so, the projection $\pi:\E \rightarrow \W$
 cannot technically be proper.}.
In~\cite{buzcal}, this was proved in the affirmative
for the particular case of tame level $N = 1$ and $p = 2$.
The proof, however, was quite explicit and required
 (at least indirectly) that the curve
$X_0(Np)$ have genus zero.
In this paper, we work with general $p$ and arbitrary tame level,
 although our result only applies
 at certain arithmetic weights in the center of weight space.

Recall that the $\C_p$-points of $\W$ are the continuous homomorphisms
from the Iwasawa algebra 
$\displaystyle{\Lambda:=\Z_p[[\lim_{\leftarrow} (\Z/Np^k\Z)^{\times}]]}$
to $\C_p$.  Let $\cyc$ denote the cyclotomic character.
Our main theorem is: 

\begin{theorem} Let $\E$ be the $p$-adic eigencurve
of tame level $N$. Let $D$ denote the closed unit disk, and
let $D^{\times}$ denote $D$ with the origin
removed. Let  $h:D^{\times} \rightarrow \E$ be a morphism 
such that $\pi \circ h$ extends to $D$.
Suppose, moreover, that
$(\pi \circ h)(0) = \kappa$, where $\kappa$ is of the form
$$\kappa = \chi^k \cdot \psi,$$
for $k \in \Z$ and $\psi$  a finite order character of conductor dividing
$N$.
Then there exists a  map $\wt{h}:D \rightarrow \E$
making the following diagram commute:
$$
\xymatrix@=4em{
     D^{\times} \ar[r]^{h} \ar[d] & \E \ar[d]^{\pi} \\
     D \ar[r] \ar[ru]_{\wt{h}} & \W \\}
      $$
\label{theorem:main}
\end{theorem}

The new idea of this paper is, roughly speaking, to specialize the family
$h$ of overconvergent modular forms to an infinitesimal neighbourhood
of the punctured point. Using the techniques
of~\cite{buzcal}, we conclude that the limiting form ``$h(0)$'' will
be an overconvergent modular form $G_0$, and thus, it suffices
to prove that the
$U_p$ eigenvalue of this form is not $0$.
However, if $U_p G_0 = 0$, 
the infinitesimal deformations of $G_0$
will have nilpotent but nonzero $U_p$ eigenvalue.
We are able
to deduce a contradiction by combining this idea with
 the philosophy of~\cite{buzcal}
 that finite
slope forms will have large radius of convergence while infinite
slope forms will have small radius of convergence.
This idea was inspired by
work of Bella\"{\i}che and Chenevier~\cite{frenchies},
who study deformations of trianguline $\Gal(\Qbar_p/\Q_p)$-representations
by considering deformations of finite dimensional
$(\varphi,\Gamma)$-modules over Artinian extensions of $\Q_p$.
Their goal was to study
the tangent spaces of eigencurves using local techniques
from $p$-adic Hodge theory.
It is plausible that the properness of the eigencurve is a global
manifestation of a purely
local theorem;  such an idea was suggested to the author
 --- at least at integral weights ---
by Mark Kisin in 2001 and was 
discussed in several tea room conversations
during the 
Durham symposium on Galois 
representations 
 in 2004 and the
eigenvarieties semester at Harvard
in 2006. However, even with current
advances in the technology of local Galois representations,
a  natural conjectural statement implying properness  has not yet
been formulated. One issue  to bear in mind is that
slightly stronger statements
one may conjecture are false. For example, there exists
a pointwise sequence of finite slope forms converging
to an infinite slope form~\cite{wascol}.

\medskip

It is a pleasure to thank Kevin Buzzard for many fruitful discussions;
the debt this paper owes to~\cite{buzcal} is clear.
I would also like to thank Matthew Emerton, Toby Gee, and Mark Kisin
for useful conversations.

\section{Overconvergent Modular Forms}
\label{section:over}

Let $N \ge 5$ be an integer co-prime to $p$;  let
 $X = X_1(N)$; and let
 $X_0(p) =
 X(\Gamma_1(N) \cap \Gamma_0(p))$.
Since $N \ge 5$, the curves
 $X$ and $X_0(p)$ are  the
compactifications
of a smooth moduli spaces. The curve $X$ comes equipped
with a natural sheaf $\omega$, which, away from the cusps,
is the pushforward of the sheaf of differentials on
the universal modular curve.
Let $A$ be a  characteristic zero lift of the Hasse invariant
with coefficients in
 $W(\Fbar_p)[[q]]$, and thus,
$A \in H^0(X/W(\Fbar_p),\omega^{\otimes (p-1)})$ by the $q$-expansion
principle.
We further insist
that $A$ has trivial character if $p > 2$,
and that $A^2$ has trivial character if $p = 2$; this
is possible since $N > 1$.
Let $X_0(p,r)  \subseteq X^{\mathrm{an}}_0(p)$ denote the connected component
containing $\infty$
of the affinoid $\{x \in X^{\mathrm{an}}_0(p); \ |A(x)| \ge |r|\}$.
 Standard arguments 
imply
that $|A(x)|$ on $X_0(p,r)$ is
independent of the choice of $A$, provided that $v(r) < p/(p+1)$.

\medskip

Let $r \in \C_p$ be an element with
 $p/(p+1) > v(r) > 0$.              
Let $\cyc$ denote the cyclotomic character;
let $\psi$ denote a finite order character of conductor dividing $N$;
and let $k \in \Z$.
\begin{df} The overconvergent modular forms of weight $\cyc^k \cdot \psi$,
level $N$, and radius of convergence $r$ are sections of
$H^0(X_0(p,r),\omega^{\otimes k})$ on which the diamond operators
act via $\psi$. We
denote this space by $M(\C_p,N,\cyc^k \cdot \psi;r)$.
The space of overconvergent modular forms
of weight $\cyc^k \cdot \psi$ and level $N$ is
$$M(\C_p,N,\cyc^k \cdot \psi) := \bigcup_{|r|<1} M(\C_p,N,\cyc^k
\cdot \psi;r).$$
\end{df}
The space $M(\C_p,N,\cyc^k \cdot \psi;r)$ has a natural
Banach space structure. If $\cyc^k  = 1$, the norm
$\|\cdot\|$ is the supremum norm.

\medskip

Let $\kappa \in \W(\C_p)$ denote a point in weight space. Recall
that the Eisenstein series $E(\kappa)$ is defined away
from zeroes of the Kubota--Leopoldt zeta function by the following
formulas:
$$E(\kappa) = 1 + \frac{2}{\zeta(\kappa)} \sum_{n=1}^{\infty}
\sigma^{*}_{\kappa}(n) q^n, \qquad
\sigma^*_{\kappa}(n) = \sum_{(d,p)=1}^{d|n} \kappa(d) d^{-1}.$$
The coefficients of $E(\kappa)$ are rigid analytic functions
on $\W$.
If $\kappa$ is trivial on the roots of unity in $\Q_p$, then,
as a $q$-expansion, $E(\kappa)$ is congruent to $1$ modulo
the maximal ideal of $\Zbar_p$.
Coleman's idea is to define overconvergent forms of
weight $\kappa$ using the formal $q$-expansion $E(\kappa)$.
Before we recall the definition, we also recall some elementary
constructions related to  weight space.
If 
$$\Z_{p,N}:= \lim_{\leftarrow} (\Z/Np^k\Z)^{\times},$$
then there is a natural isomorphism $\Z_{p,N} \simeq
(\Z/N \q \Z)^{\times} \times (1 + \q \Z_p)$, where
$\q = p$ if $p$ is odd, and $\q = 4$ otherwise.
If $a \in \Z_{p,N}$, then $\lll a \rrr$ denotes the projection
of $a$ onto the second factor, and $\tau(a) = a/\lll a \rrr$ the
projection onto the first. The rigid analytic space
$\W$ has a natural group structure. Denote the connected
component of $\W$ by $\B$; the component group of $\W$
 is $(\Z/N \q \Z)^{\times}$. If $\kappa \in \W(\C_p)$,
then let $\langle \kappa \rangle$ denote the weight
$a \mapsto \kappa(\lll a \rrr)$ and
$\tau(\kappa)$ the weight $a \mapsto \kappa(\tau(a))$; $\langle \kappa \rangle$
is the natural projection of $\kappa$ onto $\B$.
If $\cyc$ denotes the cyclotomic character, then
for any character $\psi$ of $(\Z/\q N \Z)^{\times}$, there
is a unique congruence class modulo $p-1$ (or modulo $2$ if $p = 2$)
such that for any $k \in \Z$ in this congruence class,
$\tau(\eta \cdot \cyc^{-k})$ has conductor dividing $N$.
We fix once and for all a choice of representative $k \in \Z$
for this congruence class. 

\medskip

We recall now the definition
of overconvergent modular
forms of weight $\kappa$:
\begin{df}
Overconvergent modular forms of
weight $\kappa$ and tame level $N$
 are $q$-expansions of the form
$V E_{\langle \kappa \cdot \cyc^{-k} \rangle} \cdot F$,
where 
$F \in M(\C_p,N, \cyc^k \cdot \tau(\kappa \cdot \cyc^{-k}))$.
\end{df}
Note that this is not the exact definition that
occurs on~\cite{eigencurve}, \S 2.4, since
we have chosen to work with $\Gamma_0(p)$ structure
rather than $\Gamma_1(p)$ structure.
Yet both definitions are easily
seen to be equivalent,
using, for example, Theorem 2.2.2 of \emph{ibid}.
We do not define the radius
of convergence of an overconvergent form of general weight.

\section{Hasse Invariants}

In this section, we prove some estimates for the convergence
of certain overconvergent modular forms related to Hasse invariants.
As in Section~\ref{section:over}, let $A$ be a characteristic
zero lift of the Hasse invariant with coefficients in 
 $W(\Fbar_p)[[q]]$.

\begin{lemma} Let $v(r) < 1/(p+1)$, and let $x$ be a point on
$X_0(p,r)$.
Then
$$\frac{A(x)}{\VA(x)} \equiv 1 \kern-0.5em{\mod \frac{p}{A(x)^{p+1}}}.$$
\label{lemma:one}
\end{lemma}

\begin{proof} 
The weight of $A$ is
$p -1$. Let $E$ be the elliptic curve associated
to $x$, and $H$ the canonical subgroup.
Let $\omega_E$ be a N\'{e}ron differential of
$E$, and let $a = A(E,\omega_E)$ (we implicitly
trivialize $H^0(E,\Omega^1)$).
 By Theorem 3.1 of Katz (\cite{katz}, 
p.113), we deduce that
$E/H$ is isomorphic to $E^{(p)}$ modulo $p/a$, where
$E^{(p)}$ is the image of $E$ under Frobenius. Hence,
$$A(E/H,\omega_{E/H}) \equiv a^{p} \kern-0.5em{\mod p/a},$$ 
where $\omega_{E/H}$ is any differential that
can be identified with the inverse image of $\omega_E$
under Frobenius modulo $p/a$.
By definition,
$$\VA(E,\omega_E) = p^{1-p} \cdot A(E/H,\pi^*\omega_E)
= p^{1-p} \cdot \lambda^{1-p} \cdot  A(E,\omega_{E/H}),$$
where $\pi^{*} \omega_E = \lambda \cdot \omega_{E/H}$.
Remark 3.6.5.0 of Katz (\cite{katz},  
p.116) identifies $\lambda$
with 
$$a_{p-1}/p \equiv
 A(E,\omega_E)/p \kern-0.5em{\mod 1} \equiv a/p \kern-0.5em{\mod 1}.$$
The factor of $p$ comes from the identity
 $dx^p/x^p = p(dx/x)$.
It follows directly that
$ \VA(E) \equiv a \kern-0.5em{\mod p/a^p}$, and 
the lemma follows after dividing by $a$.
\end{proof}

\medskip

\begin{cor} Suppose that $v(r) < 1/(p+1)$.  Then 
$\log(A/\VA) \in M(\C_p,N,1,r)$.  
If $s \in \C_p$ is sufficiently small, then
$(A/\VA)^s
\in M(\C_p,N,1,r)$.
\label{cor:corny}
\end{cor}

\begin{proof} From Lemma~\ref{lemma:one}, we deduce that
$A/\VA - 1$ has norm $< 1$ on $X_0(p,r)$, which implies the first claim.
Moreover,
$\|s \cdot \log(A/\VA)\| \ll 1$ for sufficiently small $s$, and hence,
if $s$ is sufficiently small,
$$(A/\VA)^s =
   \exp\left(s \cdot \log\left(A/\VA)\right)\right)$$
is well-defined and lies in 
$M(\C_p,N,1,r)$.
\end{proof}

\section{Families of Eigenforms}

Let $h: D^{\times} \rightarrow \E$ denote
an analytic family of overconvergent modular eigenforms                          
of finite slope such that $\pi \circ h$ extends                        
to $D$, and suppose that
$(\pi \circ h)(0) =\kappa$, 
 where $\kappa$ is of the form
$\kappa = \chi^k \cdot \psi$ with $k \in \Z$ and a finite
order character
$\psi$ of conductor dividing
$N$. We assume that the image of $h$ lies in
the cuspidal locus since the Eisenstein locus is easily
seen to be proper
 (cf~\cite{buzcal}, Theorem 8.2).
Any weight in $\W(\C_p)$ sufficiently close to
$\kappa$ is of the form
$\kappa \cdot \B^*$, where
$$\B^{*}:= \left\{ \eta(s): a \mapsto
\lll a \rrr^{s(p-1)} \ | \ s \in \C_p, v(s) > -1 + \frac{1}{p-1}
 \right\},$$
(the inequality should be $v(s) > - 1$ when $p = 2$).
Our definition of $\B^*$ is normalized
slightly differently from~\cite{eigencurve}~p.28,
as we have included an extra factor of $p-1$ in
the exponent. 
After shrinking $D$, if necessary, we may assume
that $(\pi \circ h)(D^{\times})
\subset \kappa \cdot \B^{*}$.
Given $t \in D$, we may consider $h(t)$ to be
a normalized
eigenform in $M(\C_p,N,\kappa \cdot \eta(s(t)))$,
for some $\eta(s(t)) \in \B^*(\C_p)$ and analytic function $s(t)$. By
assumption, $U h(t) = \lambda(t) h(t)$ for some analytic function
$\lambda(t)$ which does not vanish on $D^{\times}$.
By considering $q$-expansions, we deduce that $h(0)$
exists as a $p$-adic modular form in the sense
of Katz~\cite{katz} (for a more detailed
proof, see~\cite{buzcal}, p.229). 
The modular form $A$ has weight $\chi^{p-1} = \eta(1)$ if
$p > 2$, and $A^2$ has weight $\chi^2 = \eta(2)$ if $p = 2$.
Thus (shrinking $D$ again if necessary), we may construct
a map
$$\vpi: D^{\times} \rightarrow M(\C_p,N,\kappa)$$
via the formula
$\vpi(t) = h(t)/\VA^{s(t)}$. This map is 
well-defined as an easy consequence of Corollary B4.2.5
of~\cite{coleman}, namely that $E_s/A^s$ is overconvergent
of weight zero
where $E_s$ is the Eisenstein series of weight
$\eta(s)$.

\begin{lemma} Suppose that $v(r) < 1/(p+1)$. After
shrinking $D$, if necessary, the image of  $\vpi$ lands
in 
$M(\C_p,N,\kappa,r)$. 
\label{lemma:extend}
\end{lemma}

\begin{proof}
By construction, $\vpi(t)$ lies in
 $M(\C_p,N,\kappa,\mu)$ for some $\mu$ with $v(\mu) > 0$.
Since $\kappa$ is of the form $\chi^k \cdot \psi$, we may
therefore realize $\vpi(t)$ as a section of
$H^0(X_0(p,\mu),\omega^{\otimes k})$. Here we use the
fact that $\psi$ has conductor co-prime to $p$.
Consider the operator $\vU = (A/\VA)^{s(t)} U$,
where $U$ is the usual operator on overconvergent
modular forms~\cite{coleman,coleman1}. If $s(t)$ is sufficiently small,
then by Corollary~\ref{cor:corny},
the factor $(A/\VA)^{s(t)}$ lies in $M(\C_p,N,1,r)$. On the other hand,
$$\vU(\vpi(t)) = (A/\VA)^{s(t)} U(\vpi(t)/\VA^{s(t)}) = 
(A/\VA)^{s(t)} (\lambda(t) g(t)/A^{s(t)}) = \lambda(t)  \vpi(t).$$
If $v(\mu) < v(r)$,
then $U$ maps $M(\C_p,N,\kappa,\mu)$ to $M(\C_p,N,\kappa,\mu^{p})$.
Thus, since $\lambda(t) \ne 0$ for $t \in D^{\times}$, we deduce 
from the equality $g(t) = \lambda(t)^{-1} \vU(\vpi(t))$ that
if $\vpi(t)$ lies in $M(\C_p,N,\kappa,\mu)$, then $\vpi(t)$ 
lies in
$M(\C_p,N,\kappa,\max\{\mu^p,r\})$. Thus, by induction, $g(t)$ lies
in $M(\C_p,N,\kappa,r)$.
\end{proof}

\medskip

Let $Y$ be a connected affinoid variety,
and let $V$ be a non-empty
admissible open affinoid subdomain of $X$.
Let $B =  \Sp(\C_p\langle T\rangle)$, and $A = \Sp(\C_p\langle T,T^{-1} \rangle)$.
Let $\sheaF$ denote a sheaf on $Y$ such that $\sheaF(Y) \rightarrow \sheaF(V)$
is an inclusion.
The following is an immediate generalization of~\cite{buzcal},
Lemma 8.1.

\begin{lemma} Let $\sheaG$ be the pullback of $\sheaF$ to $Y \times B$.
If $g$ is a section of $\sheaG(V \times B)$ that extends to a section
of $\sheaG(Y \times A)$, then $g$ extends to $\sheaG(Y \times B)$.
\label{lemma:buzz}
\end{lemma}

\begin{proof} By assumption  $g \in
 \sheaG(V \times B) = \sheaF(V)\langle T \rangle$ and
$g \in \sheaG(Y \times A) = \sheaF(Y) \langle T,T^{-1} \rangle$.
Since $\sheaF(Y) \subset \sheaF(V)$, the intersection of these
two modules inside $\sheaG(V \times A) = \sheaF(V)\langle T,T^{-1} \rangle$ is
$\sheaF(Y)\langle T \rangle = \sheaG(Y \times B)$.
\end{proof}

\medskip

As remarked above, the
$q$-expansion $\vpi(0) = h(0)$ is a Katz $p$-adic
modular form of weight $\kappa$. 
Let $Y = X_0(p,r)$; 
let  $V = Y^{\mathrm{ord}}$ be the ordinary locus of $Y$;
and, let $\sheaF = \omega^{\otimes k}$.
Since $B(\C_p) = D$, the map $\vpi$
extends to a morphism $B \rightarrow \sheaF(V)$. On
the ``boundary'' $A$ of $B$ 
(or on any annulus contained in $B$ and not containing zero),
 $\vpi$ extends to a morphism
$A \rightarrow \sheaF(Y)$. Since morphisms from
$B$ to $\sheaF(Y)$
may be identified with
$\sheaG(Y \times B)$, we deduce from Lemma~\ref{lemma:buzz}
that $\vpi$ extends to a morphism $B \rightarrow \sheaF(Y)
= M(\C_p,N,\kappa,r)$.
Thus, to complete the proof of Theorem~\ref{theorem:main},
it suffices to prove
that $g(0)$ has finite slope
 or, equivalently, that 
$\lambda(0) \ne 0$. Hence, we assume that $\lambda(0) = 0$.
 Since $\lambda$ doesn't vanish
on $D^{\times}$, it is not identically zero, and thus,
$$\lambda(T) = \lambda_m T^m + \ldots,$$
for some $m \in \mathbf{N}_{>0}$ such that 
$\lambda_{m} \ne 0$.
There is, moreover, an identity
$$U_{T}(g(T)) = \exp\left(s(T) \cdot \log\left(\frac{A}{\VA}\right)
\right) U(g(T)) = \lambda(T) g(T).$$
We now specialize this identity to $\C_p[\eps]/\eps^{m+1}$
via the map $T \mapsto \eps$. This specialization  is
not strictly necessary, as one could simply
work with the first $m$ coefficients of the Taylor
expansion of $g(T)$. We persist, however,
for psychological reasons,
in order to view $g(\eps)$ as associated to a form with weight in some
infinitesimal neighbourhood of $\kappa$.
Suppose that $\displaystyle{g(\eps) = \sum_{k=0}^{m} G_k \cdot
\eps^k}$.
Then, since $s(0) = 0$, it follows that $s(\eps) \equiv 0 \mod \eps$,
and thus,
$$U g(\eps) =  \exp\left(-s(\eps) \cdot \log\left(\frac{A}{\VA}\right)
\right) \cdot \lambda(\eps) \cdot g(\eps)
= \lambda_m \eps^m  G_0.$$
By equating coefficients, we find that
$U G_0 = 0$, and $U G_m = \lambda_m G_0$.
Since $U$ increases
overconvergence (and $\lambda_m \ne 0$), it follows that
$G_0 \in M(\C_p,N,\kappa,r^p)$.
Our only condition on $r$ so far is that $v(r) < 1/(p+1)$.
Thus, we take $p \cdot v(r) = v(r^p) = 1/(p+1)$.

\begin{lemma} If $g(0) = G_0$ does not have finite
slope, then
\begin{enumerate}
\item $G_0 \in M(\C_p,N,\chi^k \cdot \psi,p^{1/(p+1)})
= H^0(X_0(p,p^{1/(p+1)}),\omega^{\otimes k})^{\langle \rangle = \psi}$.
\item $U G_0 = 0$.
\item $G_0 = q + \ldots \ne 0$.
\label{lemma:contra}
\end{enumerate}
\end{lemma}

\begin{proof} The first two claims are proved
above. For the final claim, note that $g(0)$
is a limit of normalized cuspidal eigenforms, and
so the first coefficient is $q$.
\end{proof}

To complete the proof, we note that
the conclusions of Lemma~\ref{lemma:contra}
are in contradiction with the following
result from~\cite{buzcal}.

\begin{lemma}
If $k\in\Z$ and $G\in H^0(X_0(p,p^{1/(p+1)}),\omega^{\otimes k})$
is in the kernel of $U$, then $G = 0$.
\end{lemma}

\begin{proof} Suppose $G\in H^0(X_0(p,p^{1/(p+1)})),\omega^{\otimes k})$ is
arbitrary.
Let $E$ be an elliptic curve over a finite extension
of $\Q_p$, equipped with  a subgroup $C$
of order $p$ and  with level $N$ structure $L$.
If the corresponding point $(E,C,L)\in Y$ is in 
$X_0(p,p^{1/(p+1)})$,
then one can regard $F(E,C,L)$ as an element
of $H^0(E,\Omega^1)^{\otimes k}$. Now define
$F\in H^0(X_0(p,p^{(p/(p+1))}),\omega^{\otimes k})$ by
$$F(E,L)=\sum_{D\not=C}\pi^*G(E/D,\overline{C},\overline{L}),$$
where the sum is over the subgroups~$D\not=C$ of $E$ of order~$p$;
$\pi$ denotes the projection map $E\to E/D$;
$\pi^{*}$ denotes the pullback from $H^0(E/D,\Omega^1)^{\otimes k}$
to $H^0(E,\Omega^1)^{\otimes k}$; and, a bar over a level
structure denotes its natural pushforward. An easy calculation
using Tate curves (see, for example, Proposition~5.1 of~\cite{wild})
shows that $F=p \kern+0.07em{U G}$, and hence, if $U G=0$,
 then $F=0$. 
If $E$ is an elliptic curve with no canonical subgroup, and we fix
a level $N$ structure $L$ on $E$, 
then $(E,C,L)\in X_0(p,p^{(p/(p+1))})$
 for all $C$.
Thus, $F(E,C,L)=0$
for such $E$, and hence,
$$\sum_{D\not=C}\pi^*G(E/D,E[p]/D,\overline{L})=0.$$
for all $C$. 
Summing, one deduces that
$G(E/D,E[p]/D,\overline{L})=0$ for all $D$ of order~$p$.
This implies that $G$ is identically zero on the ``boundary''
of $X_0(p,p^{1/(p+1)})$ and, hence, that $G$ is identically zero.
\end{proof}


\begin{thebibliography}{99}

\small

\bibitem{ashstevens}
A.~Ash, G.~Stevens.
$p$-adic Deformations of Cohomology on $\mathrm{GL}(n)$: the
non-ordinary case, preprint.

\bibitem{frenchies}
J.~Bella\"{\i}che, G.~Chenevier,
\emph{$p$-adic families of Galois representations
and higher rank Selmer groups}, book in preparation.

\bibitem{wild}
K.~Buzzard,
Analytic Continuation of Overconvergent Eigenforms,
J. Amer. Math. Soc. \bf 16\rm \ (2003), no.1, 29--55.

\bibitem{eigenvarieties}
K.~Buzzard,
Eigenvarieties,
To appear in the proceedings of the 2004 Durham Symposium
on $L$-functions and arithmetic.

\bibitem{buzcal}
K.~Buzzard, F.~Calegari,
The $2$-adic Eigencurve is Proper,
to appear in Documenta Mathematica, Special
Volume in Honour of John Coates.

\bibitem{chenevier}
G.~Chenevier, 
Familles $p$-adiques de formes automorphes pour ${\rm GL}\sb n$. 
J. Reine Angew. Math. 570 (2004), 143--217. 


\bibitem{coleman}
R.~Coleman,
$p$-adic Banach spaces and families of modular
forms,
Invent. Math. \bf 127\rm \ (1997) no. 3, 417--479.

\bibitem{coleman1}
R.~Coleman,
Classical and overconvergent modular forms.
Invent. Math. \bf  124\rm \  (1996), no. 1-3, 215--241.



\bibitem{wascol}
R.~Coleman, W.~Stein.
Approximation of eigenforms of infinite slope by eigenforms
of finite slope,
Geometric aspects of Dwork theory. Vol. I, II,  437--449,
2004.


\bibitem{eigencurve}
R.~Coleman, B.~Mazur.
The Eigencurve,
 Galois representations 
in arithmetic algebraic geometry (Durham, 1996),  1--113,
London Math. Soc. Lecture Note Ser., 254, 
Cambridge Univ. Press, Cambridge, 1998.

\bibitem{colmez}
P.~Colmez,
S\'{e}rie principale unitaire pour
$\mathrm{GL}_2(\Q_p)$ et
repr\'{e}sentations traingulines de dimension $2$,
Preprint.


\bibitem{emerton}
M.~Emerton,
On the interpolation of systems of eigenvalues attached to automorphic Hecke eigenforms. 
Invent. Math. \bf 164\rm \  (2006), no. 1, 1--84. 


\bibitem{katz}
N.~Katz,
$p$-adic properties of modular schemes and modular forms.  
Modular functions of one variable, III 
(Proc. Internat. Summer School, Univ. Antwerp, Antwerp, 1972),  
pp. 69--190. Lecture Notes in Mathematics, Vol. 350, Springer, Berlin, 1973. 

\bibitem{urban}
E.~Urban,
Eigenvarieties for reductive groups,
preprint.


\end{thebibliography}
\end{document}